\input amstex
\documentstyle{amsppt}
\pageheight{50.5pc} \pagewidth{32pc}

\define\1{\hbox{\rm 1}\!\! @,@,@,@,@,\hbox{\rm I}}

\topmatter
\title
{Limit raring processes with apllication}
\endtitle
\footnotetext{This research was supported (in part) by the
Ministry of Education and Science of Ukraine, project No 01.07/103
}
\author
{ Vitalii A. Gasanenko}
\endauthor

\address
{Institute of Mathematics, National Academy of Science of Ukraine,
Tereshchenkivska 3, 01601, Kiev, Ukraine}
\endaddress
\email {gs\@imath.kiev.ua}
\endemail
\subjclass {60 G 55}
\endsubjclass
\endkeywords
\subjclass {60 J 60}
\endsubjclass
\email {gs\@imath.kiev.ua}
\endemail
\abstract This paper deals with study of the sufficient condition
of approximation raring process with mixing by renewall process.
We consider use the proved results to practice problem too
\endabstract

\endtopmatter

\document

\bigskip

\bigskip

    If we have a strictly increasing  almost sure sequence  of positive
random values

$\quad\{\tau_{i},  i\geq 0\},\quad
\tau_{0} = 0 \quad\tau_{i+1} > \tau_{i}, \quad i\geq 0$ then we can
define random flow of  point-event on the time axes. The moment
appereance $i$-th event coincides with time $\tau_{i}$. Any
subflow this flow is named raring flow. Thus  $i$ - th event
in raring flow has number $\beta(i)$ in initial flow
(it is clear $\beta(i)\geq i$). We wish to investigate
the sequence $\{\beta(i), i\geq 0\}$. We make more precise
the limit theorem from [1] and consider new application too.

\bigskip

\centerline{\bf 1. Limit theorem}

\bigskip

     Let us consider the sequence of discrete random values

$$
\xi(t),\quad t\in \{0,1,2,...\},\quad \xi(t) \in \{1,2,...\}.
$$

     We are going to investigate distribution the following sequence
$$
\beta(1) = \xi(0),~~ \beta(m+1)=\beta(m)+\xi(\beta(m)),~~ m\geq 1,
$$

For this purpose, we introduce  the following objects
$$
v(t) = \max\{m\geq 1: \beta(m)\leq t\},
$$

$$
 \alpha(k) = \sup \limits_{x\geq 0}\quad\sup \limits_{A\in F_{\leq x},
B\in F_{\geq x+k}}| P(AB) - P(A)P(B) |,
$$

$$
F_{\leq x} := \sigma(\xi(s), s\leq x),\quad
F_{\geq x} := \sigma(\xi(s), s\geq x).
$$

\proclaim{\bf Statement}The following inequlity holds for any $x>0$

$$
P(\beta(m)<x)\leq \max\limits_{t\leq x}P(\xi(t)<\frac{x}{m})([x]+1).
$$
\endproclaim

\demo{\bf Proof} We have by definition of $\beta(m)$

$$
\{\beta(m)<x\}\subseteq\{ \bigcup\limits_{i=1}^{[x]+1}\{\xi(i)<\frac{x}{m}\}\},
$$

from latter one proof must.
\enddemo

    Now we will proof the limit theorem for random values $\beta(m)$ in case
when process $\xi(t)$ depends on parameter $n$.
The dependence on $n$ means , in this case, that
sequence processes $\xi_{n}(t)$ must convergence to infinity (in some sense) at fixed $t$
under $n\to \infty$. Such situation occurs in practice problem very often.

     The parameter $n$ is index for all values are defined by  $\xi_{n}(t)$ .

For example, the values  $v(t)$ transform to  $v_{n}(t)$.

Let  $\mathop{\Rightarrow}\limits_{n\to \infty}$ denotes weak convergence
of random values or disribution functions. Let $N(t)$ is equal to
number of renewals on the interval $[0,T]$ of renewal process
$\{\eta_{i}\}_{i\geq 1}$ . This process has the following property

$$
P(\eta_{1} \leq x) = R_{1}(x),\quad P(\eta_{i} \leq x)= R_{2}(x), i\geq 2.
$$

Here $R_{1}(\cdot), R_{2}(\cdot)$ are a distribution function.

\proclaim
{\bf Theorem 1}~ If sequence of numbers $c_{n} \to \infty$ under $n\to \infty$
exist such that the following conditions hold
:

1)~ $\lim\limits_{n\to \infty}P(\xi_{n}(0) c_{n}^{-1} \leq x) = R_{1}(x);$

2)~ $\lim \limits_{n\to \infty}
\sup \limits_{a\leq \delta \leq t}
\mid P(\xi_{n}([c_{n}\delta])c_{n}^{-1} \leq y) - R_{2}(y)\mid = 0;$

$\delta$ -- any positive number, ~~~~$ t< \infty$,

functions $R_{i}(y)$ are continuos functions for $y>0$;

3)~ $\lim\limits_{n \to \infty}\alpha_{n}(c_{n})c_{n} = 0$.

then

$$
v_{n}(c_{n}t) \mathop{\Rightarrow}\limits_{n\to \infty} N(t),
$$

for every fixed $t$.
\endproclaim

\demo{\bf Proof}

We denote by $\beta_{k}(m), ~~m\geq 1$ the sequence which is defined
by the  sequence  $\beta(m)$ under condition $\xi(0)=k$.

That is
$P(\beta_{k}(m)=s)=P\left(\beta(m)=s/\xi(0)=k\right)$.

Futher~~ $v(k,t)=\max\{m\geq 1: \beta_{k}(m)<t\}$.

We define the following sequence of random values $\nu_{k}(m)$:

$$
\nu_{k}(0)\equiv 0,\quad\nu_{k}(1)=\xi(k),\quad \nu_{k}(m+1)=\nu_{k}(m)+\xi(\nu_{k}(m)),
~~m\geq 1.
$$

Futher let  $V_{k}(t)=\max\{m\geq 1:~~\nu_{k}(m)\leq t \}$.

Now we introduce the sequence of random integer numbers  $\beta_{l,k}(m),~~m\geq 1$,
which have the following distribution function

$$
P(\beta_{l,k}(m)=s)=P\left(\beta(m)=s/\xi(0)=l,~\xi(l)=k\right)=
$$
$$
=P\left(\nu_{l+k}(m)=s-l-k/\xi(0)=l,~\xi(l)=k\right).
$$

We will denote by ~~$v_{l,k}(t)=\max\{m\geq 2:~~\beta_{l,k}(m)\leq t\}$.

 By the definition of  $v(t)$ and  $v(l,t)$ we have stochastic
equalities (right and left parts have the same distribution function)

$$
v(t) \doteq \sum\limits_{l=1}^{[t]}I(\xi(0) = l)(v(l,t) +1),
$$

$$
v(l,t)=\sum\limits_{k=1}^{[t-l]}I\left(\xi(l)=k/\xi(0)=l\right)(v_{l,k}(t)+1)
$$

here the function~ $I(\cdot)$ is indicator function of sets.

Applying indicator identity

$$
s^{I(\cdot) x} = 1 + I(\cdot)(s^{x} - 1),
$$

we get

$$
Ms^{v(t)} = 1 + \sum\limits_{l=1}^{[t]}MI(\xi(0) =l)(s^{v(l,t) +1} - 1),
$$

$$
Ms^{v(l,t)}=1+\sum\limits_{k=1}^{t-l}MI\left(\xi(l)=k/\xi(0)=l\right)\left(
s^{v_{l,k}(t)+1}-1\right)
$$

here $s \in (0,1).$

If the  $\xi(t)$ depends on parameter $n$, then latter equalities
have the following forms.

Put

$$
Ms^{v_{n}(c_{n}t)}=g_{n}(c_{n}t,s),\quad Ms^{V_{n,k}(c_{n}t,s)}=f_{n,k}(c_{n}t,s).
$$

Futher

$$
g_{n}(c_{n}t,s)=1-P(\xi_{n}(0)\leq t)+s\sum\limits_{l=1}^{[c_{n}t]}
MI(\xi_{n}(0)=l)
s^{v_{n}(l,c_{n}t)},
$$

$$
Ms^{v_{n}(l,c_{n}t)}=1- P\left(\xi_{n}(l)\leq c_{n}t/\xi_{n}(0)=l\right)+
$$
$$
+s\sum\limits_{k=1}^{c_{n}t-l}MI\left(\xi_{n}(l)=k/\xi_{n}(0)=l\right)
s^{v_{n,l,k}(c_{n}t)} \eqno(1)
$$

We will devide  the sums in the right parts equalities  (1) into two sums:

$$
 \sum\limits_{1}^{[c_{n}\delta]} + \sum\limits_{[c_{n}\delta] +1}^{[c_{n}t]}.
\eqno (2)
$$

\bigskip

The first sum we can make less than given numder. This follows from
the conditions 1,2 and continous of functions $R_{i}(\cdot)$ in zero.

The second sum consists of the expectations of two random factors.
These factors are bounded one and mesured with respect to $\sigma$-algebrs
$F_{\leq x}, F_{\geq x+ c_{n}\delta}$ respectively
The latter one enables us to change every summand of second part of  (2)
by factor of expectations given random values with error less than
$2\alpha_{n}(c_{n}\delta)$ (look for example (20.29)[2]):

We have the following  estimates under $l\geq c_{n}\delta,$

$$
\mid MI(\xi_{n}(0) = l)s^{v_{n}(l,c_{n}t) } -
MI(\xi_{n}(0) = l)Ms^{v_{n}(l,c_{n}t)}\mid \leq 2\alpha_{n}(c_{n}\delta),
$$

$$
Ms^{v_{n}(l,t)}=\sum\limits_{d=0}^{c_{n}t-l}s^{d}
\left(P\left(\nu_{n,l}(d)\leq c_{n}t-l,
~\nu_{n,l}(d+1)>c_{n}t-l/\xi_{n}(0)=l\right)\pm\right.
$$

$$
\left.\pm P(\nu_{n,l}(d)\leq c_{n}t-l,~\nu_{n,l}(d+1)>c_{n}t-l)\right)=
Ms^{V_{n,l}(c_{n}-l)}+\pi_{n}.
$$

here ~~ $|\pi_{n}|\leq K \alpha_{n}(c_{n}\delta),~~K<\infty$.

$$
\left|P\left(\xi_{n}(l)\leq c_{n}t/\xi_{n}(0)=l\right)-
P(\xi_{n}(l)\leq c_{n}t)\right|
\leq \alpha_{n}(c_{n}\delta).
$$

Futher we have estimates in case when $k\geq c_{n}\delta$

$$
\left|MI\left(\xi_{n}(l)=k/\xi_{n}(0)=l\right)s^{v_{n,l,k}(c_{n}t)}-\right.
$$
$$
\left.-
MI\left(\xi_{n}(l)=k/\xi_{n}(0)=l\right)
Ms^{v_{n,l,k}(c_{n}t)}\right|
\leq 2 \alpha_{n}(c_{n}\delta);
$$

$$
\left|Ms^{v_{n,l,k}(c_{n}t)} - Ms^{V_{n,l+k}(c_{n}t-l-k)}\right|\leq K_{1}\alpha_{n}(c_{n}
\delta),~~ K_{1}<\infty.
$$

Now we can rewrite  (1) in the following form
$$
g_{n}(c_{n}t,s) = 1 - P(\xi_{n}(0)\leq c_{n}t) + a_{n,1}(\delta) + b_{n,1}+
s\sum\limits_{l=[c_{n}\delta]+1}^{[c_{n}t]}P(\xi_{n}(0) = l)f_{n,l}(c_{n}t-l),
$$

$$
f_{n,l}(c_{n}t-l ) = 1 - P(\xi_{n}(l)\leq c_{n}t) + a_{n,2}(\delta) +
b_{n,2} +
$$

$$
+ s\sum\limits_{k=[c_{n}\delta]+1}^{[c_{n}t]}P(\xi_{n}(l)=k)
f_{n,l+k}(c_{n}t -l-k),\quad l\geq [c_{n}\delta], eqno(3.35)
$$

here

$$
\mid b_{n,i}\mid\leq k_{i}\alpha_{n}(c_{n}\delta)[c_{n}t] ,~~k_{i}<\infty~;
\quad
 a_{n,1}(\delta) \leq P(0<\xi_{n}(0) \leq c_{n}\delta),
$$

$$
 a_{n,2}(\delta) \leq \sup\limits_{q\geq [c_{n}\delta]}
P(0<\xi_{n}(q)\leq c_{n}\delta).
$$

Futher we introduce a sequence of independence random values with the same
distribution function  ~~$\{\eta_{k}(n,\delta)\}_{k\geq 1}$~~under fixed
$\delta$. The distribution function is difined the following equlity

$$
P(\eta_{1}(n,\delta)\leq x)=P(\xi_{n}(c_{n}\delta)\leq x).
$$

We will denote

$$
S_{m}(n,\delta)=\sum\limits_{k=1}^{m}\eta_{k}(n,\delta),\quad
D_{n,\delta}(t)
=\sup\limits_{m\geq 1}\{m:~~S_{m}(n,\delta)\leq t\}.
$$
$$
Ms^{D_{n,\delta}(t)}=\sum\limits_{d=0}^{\infty}s^{d}P(D_{n,\delta}(t)=d)=:
F_{n,\delta}(t,s).
$$

We will estimate of difference of
 ~$f_{n,l}(c_{n}t-l,s),~~ l\geq c_{n}\delta$~~ and
$F_{n,\delta}(c_{n}t-l,s)$.

The definition leads to

$$
f_{n,l}(c_{n}t-l,s)=\sum\limits_{d=0}^{[c_{n}t-l]}s^{d}P(\nu_{n,l}(d)\leq
c_{n}t-l, ~\nu_{n,l}(d+1)>c_{n}t-l).
$$

Futher we get  for  ~$d=0$ by assumption  2

$$
P(\xi_{n}(l)> c_{n}t-l)\pm P(\xi_{n}(c_{n}\delta)>c_{n}t-l)=
$$
$$
=\theta_{n}+P(\xi_{n}(c_{n}\delta)>c_{n}\delta - l).
$$

Late on the designation ~$\theta_{n}$~ means that we have some
sequence of number sush that it convergence to zero under ~$n\to\infty$
and the following condition holds

$$
|\theta_{n}|\leq 2
\sup\limits_{y\leq t}\sup \limits_{\delta\leq \triangle \leq t}
\mid P(\xi_{n}([c_{n}\triangle])c_{n}^{-1} < y) - R_{2}(y)\mid
$$

We have for  ~$d=1$:~~~ $P(\nu_{n,l}(1)\leq c_{n}t-l,~\nu_{n,l}(2)>c_{n}t-l)=$
$$
=\sum\limits_{k=1}^{[c_{n}t-l]}P(\xi_{n}(l)=k,~\xi_{n}(k+l)>c_{n}t-l-k)=
$$
$$
=a_{n,\delta}+ r_{n,1,\delta}+
\sum\limits_{k=[c_{n}\delta]}^{[c_{n}t-l]}P(\xi_{n}(l)=k)
P(\xi_{n}(k+l)>c_{n}t-l-k)=
$$
$$
=a_{n,\delta}+ r_{n,1,\delta} +\theta_{n}+
\sum\limits_{k=[c_{n}\delta]}^{[c_{n}t-l]}(P(\xi_{n}(l)=k)
P(\eta_{2}(n,\delta)>c_{n}t-l-k)=
$$
$$
=a_{n,\delta}+ r_{n,1,\delta} +\theta_{n}+ P(D_{n,\delta}(c_{n}t-l)=1)-
$$
$$
- \sum\limits_{k=[c_{n}\delta]}^{[c_{n}t-l]}
\sum\limits_{s=[c_{n}\delta}^{k}(P(\xi_{n}(l)=s)-P(\eta_{1}(n,\delta)=s))
P(\eta_{2}(n,\delta)=c_{n}t-l-k)=
$$
$$
=P(D_{n,\delta}(c_{n}t-l)=1)+\pi_{n,1}.
\eqno(3)
$$

Here

\noindent $|\pi_{n,1}|\leq 2(a_{n,\delta}+\alpha_{n}(c_{n}\delta)+\theta_{n}),
~~ a_{n,\delta}=\max\limits_{i}(a_{n,i}(\delta)),~~|r_{1,n}|
\leq 2\alpha_{n}(c_{n}\delta).$

We used Abel trasformation for sum of pair factor of (3).
Similar considerations apply to case ~$d=2$~. Thus applying  (3) we get

$$
P(\nu_{n,l}(2)\leq c_{n}t-l,~\nu_{n,l}(3)>c_{n}t-l)=
a_{n,\delta}+ r_{n,2}+
$$
$$
+\sum\limits_{k=[c_{n}\delta]}^{[c_{n}t-l]}
P(\xi_{n}(l)=k)P(\nu_{n,l+k}(1)\leq c_{n}t-l-k,~\nu_{n,l+k}(2)>c_{n}t-l-k)=
$$
$$
=
P(D_{n,\delta}(c_{n}t-l)=2)-
\sum\limits_{k=[c_{n}\delta]}^{[c_{n}t-l]}
\sum\limits_{s=[c_{n}\delta}^{k}(P(\xi_{n}(l)=s)-P(\eta_{1}(n,\delta)=s))
\times
$$
$$
\times
(P(D_{n,\delta}(c_{n}t-l-k-1)=1)-
P(D_{n,\delta}(c_{n}t-l-k)=1)))
+a_{n,\delta}+ r_{n,2}+\pi_{n,1}=
$$
$$
= P(D_{n,\delta}(c_{n}t-l)=2)+\pi_{n,2}.
$$

For the latter one we used Abel transform and following
equility which checks easy.

$$
P(D_{n,\delta}(c_{n}t-l-k-1)=1)-
P(D_{n,\delta}(c_{n}t-l-k)=1))=
$$
$$
=P(S_{2}(n,\delta)=[c_{n}t]-k-l)-P(\eta_{1}(
n,\delta)=[c_{n}t]-k-l)
$$

The implicit introduced sequences have obvious sense and
the following estimates take place
 $|r_{n,2}|\leq 2\alpha_{n}(c_{n}\delta)$,
$|\pi_{n,2}| \leq 4(a_{n,\delta}+
\alpha_{n}(c_{n}\delta)+\theta_{n})$.

It is no difficult to show with help induction that we have
for ~$d=p$~  the following formulas

$$
P(\nu_{n,l}(p)\leq c_{n}t-l,~\nu_{n,l}(p+1)>c_{n}t-l)=
$$

$$
=P(D_{n,\delta}(c_{n}t-l)=p)+\pi_{n,p};\quad
|\pi_{n,p}| \leq 2p(a_{n,\delta}+\alpha_{n}(c_{n}\delta)+\theta_{n}).
$$

Thus we obtained for fixed ~$s\in (0,1)$~ representations

$$
f_{n,l}(c_{n}t-l,s)= F_{n,\delta}(c_{n}t-l,s)+L_{n},\quad
|L_{n}|\leq L(a_{n,\delta}+\alpha_{n}(c_{n}\delta)+o_{n}(1)),~~L<\infty.
$$

$$
g_{n}(c_{n}t,s) = 1 - P(\xi_{n}(0)\leq c_{n}t) + L_{n,\delta}+
s\sum\limits_{l=[c_{n}\delta]+1}^{[c_{n}t]}P(\xi_{n}(0) = l)F_{n,\delta}
(c_{n}t-l,s),
$$

$$
F_{n,\delta}(c_{n}t-l,s) = 1 - P(\xi_{n}(c_{n}\delta)\leq c_{n}t) +
Z_{n,\delta}+
$$

$$
+ s\sum\limits_{k=[c_{n}\delta]+1}^{[c_{n}t-l]}P(\xi_{n}(c_{n}\delta)=k)
F_{n,\delta}(c_{n}t -l-k),\quad l\geq [c_{n}\delta],
$$

The constructions of ~$L_{n,\delta}$~ and ~$Z_{n,\delta}$ now leads to

$$
\lim\limits_{n\to\infty}L_{n,\delta}=l_{1}L_{\delta};\quad
\lim\limits_{n\to\infty}Z_{n,\delta}=l_{2}Z_{\delta};\quad \max(l_{1},l_{2})
<\infty.
$$

here ~~$|L_{\delta}|\leq 2(R_{1}(\delta)-R_{1}(0)),\quad
|Z_{n,\delta}|\leq 2(R_{2}(\delta)-R_{2}(0)).$

Combining construction of ~$F_{n,\delta}(t,s)$~ and condition 2
we conclude that the following limit exists

$$
\lim\limits_{\delta\to 0}\lim\limits_{n\to\infty}F_{n,\delta}(c_{n}t,s)=
F(t,s),
$$

This limit is solution the following equation

$$
F(t,s)=1-R_{2}(t)+sR_{2}(\cdot)\ast F(t,s).\eqno(4)
$$

The sequence of generating function ~$g_{n}(c_{n}t,s)$~ has
limit too

$$
\lim\limits_{\delta\to 0}\lim\limits_{n\to\infty}g_{n}(c_{n}t,s)=g(t,s),
$$

This limit is solution of the following equation

$$
g(t,s)=1-R_{1}(t)+s R_{2}(\cdot)\ast F(t,s).
$$

The latter one and  (4) lead to proof of theorem.
\enddemo

\bigskip

\bigskip

{\bf Remark 1.}~ We consider a extension of theorem 1.
It consists in definition more weakly the mixing coefficient
than $\alpha(c_{n})$.

Suppose that sequence $c_{n}, n\geq 1$ from theorem1 is defined.
Now we take any sequence  $r_{n}, n\geq 1$ which satisfy
the folowing condition
$r_{n}\to \infty, r_{n} = o(c_{n})$ under $n\to \infty$ .

Futher we construct trancated process:

$$
\bar\xi_{n}(t) =
\cases
\xi_{n}(t),& \xi_{n}(t) \leq c_{n} - r_{n},\\
c_{n} - r_{n},& \xi_{n}(t) > c_{n} - r_{n}.
\endcases
$$

and construct $\sigma$- algebra
 $F_{\leq x}^{r_{n}} = \sigma(\bar \xi_{n}(t), t\leq x)$ too.

Now we define new mixing coefficient

$$
\alpha_{\xi_{n},r_{n}} = \sup\limits_{x\geq 0} \sup\limits_{ A\in F_{\leq x}
^{r_{n}}, B\in F_{\geq x + c_{n}}}\mid P(AB) - P(A)P(B)\mid .
$$

Thus this coefficient is constructed only on those events from
  $F_{\leq x}$ on which the process  $\xi_{n}(t)$ under $t\leq x$
less than value $c_{n}-r_{n}$.
Such coefficient is useful in those cases when time dependence  is controled
by values of process $\xi_{n}(t)$. For example, the event
$\{\xi_{n}(x) = k\}$ determine the behavior of process no far  interval
 $[0, x+k].$

Now we devide second sum of  (2)  in this way:

$$
\sum\limits_{[c_{n}\delta] + 1}^{[c_{n}t]} =
\sum\limits_{[c_{n}\delta] + 1}^{[(c_{n} - r_{n})t]} +
\sum\limits_{[(c_{n} - r_{n})t]+1}^{[c_{n}t]}.\eqno (5)
$$

\bigskip

We can do second sum from  (5)  less any given value due to
 continuity of $R_{i}(\cdot)$.

Futher we apply the transformation from theorem 1 to first sum with
use coefficient  $\alpha_{\xi_{n},r_{n}}(c_{n})$.

Thus we can replace the  condition 3 of theorem 1 the following condition

3'~ it exists such sequence $r_{n}: r_{n}\to \infty, r_{n} = o(c_{n})$ that

$$
c_{n}\alpha_{\xi_{n},r_{n}}(c_{n})
\mathop{\longrightarrow}\limits_{n\to \infty} 0.
$$

\bigskip

{\bf Remark 2.}If the process  $\{\tau_{i},\quad i\geq 0\}$ be such
that

$$
\lim\limits_{i\to\infty}i^{-1}\tau_{i}=
 \mu^{-1} < \infty , \quad \mu=\hbox{const}
$$

then we get convergence under conditions theorem 1

$$
P(\tau_{\beta_{n}(i)} \leq x c_{n}) \mathop{\Rightarrow}\limits_{n\to \infty}
R_{1}\ast R_{2}^{\ast (i-1)}(x \mu).
$$
It follows from known theorem of transfer (look, for example [3])

\bigskip

\centerline{\bf 2. Interaction of two renewall processes.}

\bigskip

The model of raring process which is concidered below is result interaction two
renewall process. This model was offered in  [4] as the matematical model
of practice problem.

    Let us denote by  $Z$ and $H$ two renewall processes :
$Z=\{\zeta_{i}, i\geq 1\},\quad H = \{\eta_{i}, i\geq 1 \}$.

We define stochactic characteristics of $H,\, Z$

$$
\tau_{i}=\sum\limits_{l=1}^{i}\eta_{l},\quad \vartheta_{i}=\sum\limits_{l=1}^{i}
\zeta_{l},\quad i=1,2,\dots,
$$

$$
N_{1}(t)=\sup\{n:\tau_{n}<t\},\quad N_{2}(t)=\sup\{n:\vartheta_{n}<t\},
$$

$$
\gamma_{1}^{+}(t)=\tau_{N_{1}(t)+1}-t,\quad
\gamma_{2}^{+}(t)=\vartheta_{N_{2}(t)+1}-t,\quad t>0.
$$

The points $\tau_{i}, \vartheta_{i}, i\geq 1$ are named renewall points
processes $H$ and $Z$ respectively.

If we have a renewall points of  process $Z$ in  interval
 $(\tau_{n-1},\tau_{n}]$ then
we will say that the renewall point   $\tau_{n}$ is marked by
process  $Z$.
The process $H$ marks a points of renewall of process $Z$ analogy.

Let us denote by $T_{0}^{\prime\prime}=0,T_{1}^{\prime\prime},\dots$
renewall point of  $H$ marked by  $Z$ and
$T_{1}^{\prime},T_{2}^{\prime},\dots$
renewall points of $Z$ marked by $H$.
It is clear that the following equalities
$T_{0}^{\prime\prime}=0<T_{1}^{\prime}\leq T_{1}^{\prime\prime}\leq
T_{2}^{\prime}\leq\dots$ take place.
It is shown in [4] that sequence
random values

$$
V_{n}=T_{n}^{\prime}-T_{n-1}^{\prime\prime},\quad
U_{n}=T_{n}^{\prime\prime}-T_{n}^{\prime},\quad n=1,2,\dots
$$

be Markov chain.
This chain is defined by transition probabilities

$$
P(V_{1}<x)=P(\zeta_{1}<x),
$$

$$
P(U_{n}<x|V_{n}=y),\quad P(V_{n+1}<x|U_{n}=y), \quad n= 1,2, \dots
$$

It is easy to see that for investigation
 $V_{n},\, U_{n}$ it is necessary to observe two raring processes
simultaneously:

$$
 T^{\prime\prime}=
\{T_{0}^{\prime\prime}=0,T_{1}^{\prime\prime}, T_{2}^{\prime\prime}-
T_{1}^{\prime\prime}\dots\},\quad
 T^{\prime}=
\{T_{1}^{\prime},T_{2}^{\prime}-T_{1}^{\prime},\dots\}
$$

We will investigate these raring processes separately.
We will use that the proceses  $ T^{\prime\prime}\,,\, T^{\prime}$
are raring processes respect to processes  $H\,,\,Z$ respectively.

We take, for example,  $T^{\prime\prime}$. It is define
$T^{\prime\prime}$ as subflow of  $H$ by the following indicators

$$
\chi(i) =\cases
1,&\hbox{if  $i$-th renewall point $H$  belong to}\quad  T^{\prime\prime},\\
0& \hbox{otherwise}.
\endcases
$$

$$
\xi(l) = \min\{j\geq 1: \chi(l+j) = 1\}, ~~ l\geq 0.
$$

Thus $\beta(i) = \beta(i-1) + \xi(\beta(i-1)), ~i\geq 1$ be number $i$-th
 event from  $H$ which belongs to $T^{\prime\prime}$. The moment
 $\tau_{\beta(i)}$ is moment of appearance this event.

We shall suposse that processes $H$ and $Z$ depend on a parameter
 $n,\, n\to\infty$ such that
$H_{n}=\{\eta_{n,i},\,i\geq 1\}$\quad $Z_{n}=\{\zeta_{n,i},\, i\geq 1\}$.
Now the characterictics  these processes have forms:

 $\tau_{n,i},\quad \vartheta_{n,i},\, i\geq 1,\quad
\gamma_{n,k}^{+}(t),\, t\geq 0,\, k=1,2.$

\proclaim{\bf Theorem 2} If the following conditions:

1)~there are a positive numbers $c_{n}\to \infty$ and distribution function
$G(x),\, x\geq 0$ guaranteeing the following limit

$$
\lim\limits_{n\to\infty}\sup\limits_{t\geq 0}|P(\gamma_{n,2}^{+}(t)<
\tau_{n,[nx]})\, - \, G(x)|=0,
$$

\bigskip

here $x$ is point of continuous of $G(x)$;

\bigskip

2)~ $\lim\limits_{n\to\infty} c_{n}^{-1}\tau_{n,c_{n}}
=\mu\quad \mu=\hbox{const}$.

hold then

$$
P(\tau_{\beta_{n}(k)}<x c_{n})\mathop{\longrightarrow}\limits_{n\to\infty}
G^{{\ast}k}(\frac{x}{\mu}).
$$
\endproclaim

\demo{\bf Proof}~ We will check all conditions of
theorem 1 for process $\xi_{n}(l)$.
We calculate  probability $P(\xi(l)=m)$.

$$
P(\xi(l)=1)=P(\gamma^{+}_{2}(\tau_{l})<\eta_{l+1}).
$$

$$
P(\xi(l)=2)=P(\gamma^{+}_{2}(\tau_{l})\geq \eta_{l+1},\, \gamma_{2}^{+}
<\eta_{l+2})=
$$

$$
=P(\gamma_{2}^{+}(\tau_{l})\geq \eta_{l+1},\, \gamma_{2}^{+}(\tau_{l})-
\eta_{l+1}<\eta_{l+2})=
$$

$$
= P(\gamma_{2}^{+}(\tau_{l})<\eta_{l+1}+\eta_{l+2})-
P(\gamma_{2}^{+}(\tau_{l})<\eta_{l+1}).
$$

$$
\cdots
$$

$$
P(\xi(l)=k)=P(\gamma_{2}^{+}(\tau_{l})<\eta_{l+1}+ \cdots +
\eta_{l+k})\,-\, P(\gamma^{+}_{2}(\tau_{l})<\eta_{l+1}+ \cdots +
\eta_{l+k-1}).
$$

Thus

$$
P(\xi(l)\leq m)=\sum\limits_{k=1}^{m}P(\xi(l)=k)=
P(\gamma_{2}^{+}(\tau_{l})<\eta_{l+m}+ \cdots +\eta_{l+1}).
\eqno (6)
$$

The latter one and condition 1 lead to the following
convergence

$$
P(\xi_{n}(l)< c_{n} x)=\int\limits_{0}^{\infty}
P(\gamma_{2}^{+}(t)<\tau_{n,[c_{n}x]})P(\tau_{n,l}\in dt)
\mathop{\longrightarrow}\limits_{n\to\infty}
$$

$$
\mathop{\longrightarrow}\limits_{n\to\infty}
G(x),\quad l=0,1,2,\dots.
$$

here $x>0$ is point of continuous of $G(x)$.

We have the following  equality when it is considered that (6) holds

$$
P(\xi(l)\leq m,\, \xi(l+r)\leq s) -
P(\xi(l)\leq m)P(\xi(l+r)\leq s)= 0, \, \hbox{if}\, m<r.
$$

\noindent Now we have for any sequences of numbers
$r_{n}$ such that $r_{n}\to\infty,\, r_{n} < c_{n}, n \geq 0$

$$
\alpha_{\xi_{n},r_{n}}(c_{n})= 0, \quad n\geq 0
$$

Thus all conditions of theorem 1 hold respect to process $\xi_{n}(l)$.
Now the statement of theorem 2 becomes apparent if it is remembered
the theorem of transfer.
\enddemo

Now we consider example of definition of sequence $c_{n}$
and bound function $G(x)$.

We shall suppose  that process $Z$ is Poisson process with parameter
 $\lambda_{n}$ such that $\lambda_{n}\to 0$~ under~ $n\to\infty$.
The process  $H$ don't depens on parameter $n$. It has
the following expactation of renewall interval
$\mu=M\eta_{1}<\infty$.

All these suppositions led to formula

$$
P(\gamma^{+}_{n,2}(\tau_{l})<\eta_{l+1}+\cdots +\eta_{l+m})=
\int\limits_{0}^{\infty}\lambda_{n}e^{-\lambda_{n}y}P(\tau_{m}>y)dy
=:G_{n}(m).
$$

If we put  $m=[c_{n}x]$ and make change of variables
 $\lambda_{n}y=z$  then we get

$$
G_{n}([c_{n}x])=\int\limits_{0}^{\infty}e^{-z}P(\lambda_{n}\tau_{[c_{n}x]}>z)
dz.
$$

Put $c_{n}:=\lambda^{-1}_{n}$. The indicator of set $A$ will be denoted
by $I(A)$. The following convergences are based on law of large numbers.

$$
G_{n}([c_{n}x])=\int\limits_{0}^{\infty}e^{-z}
P\left(x\frac{\tau_{[\lambda_{n}^{-1}x]}}{x\lambda^{-1}_{n}}>z\right)
\mathop{\longrightarrow}\limits_{n\to\infty}
$$

$$
\mathop{\longrightarrow}\limits_{n\to\infty}
\quad\int\limits_{0}^{\infty}
e^{-z}I(x\mu>z)dz=1-e^{-\mu x},\quad x\geq 0,
$$

and

$$
\lim\limits_{n\to\infty}\frac{\tau_{n}}{n}=
\mu.
$$

Thus all conditions of theorem 2 was cheked. The function
$G(x)$ (from condition 1 of theorem  2) be bound for
the function $G_{n}([x\lambda_{n}^{-1}])$.
In this example the moment of appereance $k$-th event in flow
{\bf $T^{\prime\prime}_{n}$} has the following function
distribution

$$
P\left(\tau_{\beta_{n}(k)}<x\lambda_{n}^{-1}\right)
\mathop{\longrightarrow}\limits_{n\to\infty}
\Bigl(1-\exp(\cdot)\Bigr)^{\ast k}(x),\quad x\geq 0.
$$

It is clear that similar example we may concider for
process $T^{\prime}$. In this case the process
 $H$  must be Poisson with "rare" events and  the process $Z$
must be a simple renewall process with bounded expectation of
time between neighboring  renewall point.


\Refs \ref \no{1} \by V.A. Gasanenko  \paper
limit theorem for raring processes with mixing I.\jour Ukrainian
Mathematical Journal \vol 50 \yr 1998 \page 471-475
\endref
\ref \no{2} \by P. Billingsley \paper Convergence of probability
measures \publ Nauka \publaddr Moscow \yr 1977 \page 352 p
\endref
\ref \no{3} \by A.A. Anisimov \paper Random processes with
discrete component\publ Vyshaij shkola \publaddr Kiev \yr 1988 \page 184
p
\endref
\ref \no{4} \by I. Kopocinska \paper Two mutually rarefied
renewall processes \jour Application Mathtmatice \vol 22 \yr
1994 \page 267-273
\endref
\endRefs
\end{document}